\documentclass[a4paper, 11pt, reqno]{amsart}
\usepackage{amssymb}
\usepackage{a4wide}
\usepackage[english]{babel}

\renewcommand{\le}{\leqslant}
\renewcommand{\ge}{\geqslant}
\newcommand{\bad}{\mathbf{Bad}}

\newcommand{\RR}{\mathbb{R}}
\newcommand{\ZZ}{\mathbb{Z}}
\newcommand{\QQ}{\mathbb{Q}}

\newcommand{\NN}{\mathbb{N}}

\newcommand{\CCC}{\mathcal{C}}
\newcommand{\DDD}{\mathcal{D}}

\newtheorem{theorem}{Theorem}

\newtheorem*{SCH1}{Theorem S1}
\newtheorem*{SCH2}{Theorem S2}

\newtheorem*{MLconj}{Mixed Littlewood Conjecture}
\newtheorem*{MSconj}{Mixed Schmidt Conjecture}
\newtheorem*{PMSconj}{Conjecture 1}
\newtheorem*{MBPV}{Conjecture 2}
\newtheorem*{BPV}{Theorem BPV}

\usepackage{amsrefs}

\begin{document}
\title{The mixed Schmidt conjecture in the theory of Diophantine approximation}

\author{Dzmitry Badziahin, Jason Levesley}

\thanks{DB: Research supported by EPSRC grant EP/E061613/1}

\author{Sanju Velani}

\thanks{SV: Research supported by EPSRC grants EP/E061613/1 and
EP/F027028/1 }

\begin{abstract}
Let $\mathcal{D}=(d_n)_{n=1}^\infty$ be a bounded sequence of integers with $d_n\ge 2$ and let $(i, j)$ be a pair of strictly positive numbers with $i+j=1$.  We prove that the set of
$x \in \RR$  for which there exists some constant $c(x) > 0$ such that
\[
\max\{|q|_\DDD^{1/i}, \|qx\|^{1/j}\}  \, >  \,   c(x)/ q    \qquad  \forall \  \  q \in \NN
\]
is one quarter winning (in the sense of Schmidt games). Thus the intersection  of any countable number
of such sets  %with differing $(i,j)$ pairs
is of full dimension. In turn, this establishes  the natural analogue of Schmidt's conjecture  within the framework of the  de
Mathan-Teuli\'{e} conjecture -- also known as the `Mixed Littlewood
Conjecture'.
  \\[0ex]


%\vspace*{1ex}
%
%{\it Keywords and phrases:} Diophantine approximation, mixed Littlewood conjecture,\newline
%\phantom{\it Keywords and phrases: ss} ~Schmidt games.
\end{abstract}

\maketitle

\section{Introduction}
\label{intro}
The famous Littlewood conjecture in the theory of simultaneous
Diophantine approximation  dates back to the 1930's and  asserts
that for every $(x,y)  \in \RR^2$, we have that
\begin{equation} \label{little}
\liminf_{q\to\infty}q\|qx\|\|qy\|=0.
\end{equation}
Here and throughout, $\| \, . \, \| $      denotes the distance to
the nearest integer.    Despite concerted efforts over the
years the conjecture remains open.   For background and recent `progress'   concerning  this fundamental  problem  see \cite{ELK,PVL} and references within.

The Schmidt conjecture in the theory of simultaneous
Diophantine approximation  dates back to the 1980's and is linked to  Littlewood's conjecture.  Given a pair
of real numbers $i$ and $j$ such that
\begin{equation}\label{neq1}
0< i,j< 1   \quad {\rm and  } \quad  i+j=1  \, ,
\end{equation}
let  $\bad(i,j)$ denote the badly approximable  set
of   $(x,y) \in \RR^2$ for which  there exists a  constant $c(x,y)> 0$ such that
\begin{equation*}
  \max \{ \; \|qx\|^{1/i} \; , \ \|qy\|^{1/j} \,  \} \, > \,
  c(x,y) \ q^{-1} \quad \forall \  q \in \NN   \ .
\end{equation*}

\noindent
A consequence of the main result in  \cite{BPV} is the following statement. Throughout,   $\dim  X $ will denote the Hausdorff dimension  of the set $X$.

\begin{BPV}
Let  $(i_t,j_t)$ be a countable number of pairs of real numbers
satisfying \eqref{neq1}.
Suppose that $ \liminf_{t \to \infty }  \min\{ i_t,j_t  \}  > 0   $.
Then
$$
\dim \Big(\bigcap_{t=1}^{\infty} \bad(i_t,j_t) \Big) = 2   \ .
$$
\end{BPV}

\noindent Thus, the  intersection of any finitely  many badly
approximable sets $\bad(i,j)$ is  trivially non-empty and therefore
establishes the following conjecture of Wolfgang M. Schmidt \cite{schconj}. {\em For any $(i_1,j_1)$ and
$(i_2,j_2)$ satisfying \eqref{neq1}, we have that}
 $$ \bad(i_1,j_1) \cap \bad(i_2,j_2) \  \neq \  \emptyset \ .
 $$

\noindent To be precise, Schmidt stated the specific problem with
$i_1=j_2=1/3$ and $j_1=i_2=2/3$. As noted by Schmidt, a counterexample to his  conjecture would  imply Littlewood's conjecture. Indeed,
the same conclusion is valid if there exists any  countable collection
of  pairs $(i_t,j_t)$ satisfying \eqref{neq1} for which the intersection of the sets $ \bad(i_t,j_t) $ is empty.

Recently, de Mathan and  Teuli\'{e} in  \cite{MT04} proposed the following variant of Littlewood's conjecture.  Let
$\mathcal{D} $ be a bounded sequence $(d_n)_{n=1}^\infty$ of
integers greater than or equal to $2$ and  let
\[
 D_0 := 1     \quad {\rm and}  \quad    D_n := \prod_{k=1}^n d_k   \ .
\]
 Now set
\[
    \omega_\mathcal{D}: \mathbb{N}\to\mathbb{N}: q \mapsto \sup\{n\in\mathbb{N}: q \in D_n\mathbb{Z} \}
\]
and
\[
    |q|_\mathcal{D} := 1/D_{\omega_\mathcal{D}(q)} = \inf\{ 1/D_n : q \in D_n\mathbb{Z} \}.
\]
When $\mathcal{D}$ is the constant sequence equal to a  prime number $p$,  the norm $|\, \cdot \, |_\mathcal{D}$ is the usual $p$-adic norm.   In analogy with Littlewood's conjecture we have the following statement.
\begin{MLconj}
 For every real number $x$
\[
     \liminf_{q\to\infty}q \, |q|_\mathcal{D} \, \|qx\|   =  0  \, .
\]
\end{MLconj}

\noindent As with the classical Littlewood conjecture, this attractive problem  remains open. The current state of affairs regarding the mixed  conjecture  is very much comparable to that of the classical one.
For background and results related to the mixed Littlewood conjecture see \cite{BDM,BHV,BM,EK,HaTa}.

It is somewhat surprising  that the analogue of  Schmidt's
conjecture within the  `mixed' framework  has to date escaped
attention.   The goal of this paper is to investigate such a
problem.   Given  $\mathcal{D} $ as above and a pair of real numbers
$i$ and $j$ satisfying \eqref{neq1}, let
\begin{equation}\label{def_badd}
\bad_\DDD(i,j):=\left\{x\in \RR :  \exists \; c(x)>0  \ { \rm so \ that \ } \max\{|q|_\DDD^{1/i}, ||qx||^{1/j}\}>  \frac{c(x)}{q}   \quad  \forall \, q\in
\NN \right\}  \,  .
\end{equation}

A consequence of the Khintchine-type result established in \cite{HaTa} is that  $\bad_\DDD(i,j)$ is of Lebesgue measure zero. The following represents a natural analogue of Schmidt's conjecture.

\begin{MSconj}
 For any $(i_1,j_1)$ and
$(i_2,j_2)$ satisfying \eqref{neq1}, we have that
 $$ \bad_\DDD(i_1,j_1) \cap \bad_\DDD(i_2,j_2) \  \neq \  \emptyset \ .
 $$
\end{MSconj}

%\begin{Schconj} For any $(i_1,j_1)$ and
%$(i_2,j_2)$ satisfying \eqref{neq1}, we have that
% $$ \bad(i_1,j_1) \cap \bad(i_2,j_2) \  \neq \  \emptyset \ .
% $$
% \end{Schconj}
%
\noindent It is easily seen that a counterexample to this conjecture would  imply the mixed Littlewood conjecture. Indeed,
the same conclusion is valid if there exists any  countable collection
of  pairs $(i_t,j_t)$ satisfying \eqref{neq1} for which the intersection of the sets $ \bad(i_t,j_t) $ is empty.  The following  constitutes our main result.

\begin{theorem}\label{SVmixschdt}
For any $(i,j)$  satisfying \eqref{neq1}, the set $\bad_\DDD(i,j)$ is $1/4$-winning.
\end{theorem}

\noindent

A consequence of {\rm winning} is the following full dimension result which settles the mixed Schmidt conjecture. See \S\ref{games} below  for the definition  and relevant implications  of winning sets.

\begin{theorem}\label{mixschdt}
For each  $t \in \NN$,  let $\mathcal{D}_t $ be a bounded sequence
as above and $(i_t,j_t)$ be a sequence of pairs of real numbers
satisfying \eqref{neq1}. Then
\[
    \dim\Big( \bigcap_{t=1}^\infty{}\bad_{\DDD_t}(i_t,j_t) \Big) = 1  \ .
\]
\end{theorem}

\noindent In a nutshell,  we are able to establish the mixed analogue of Theorem~BPV  without the annoying `$\liminf$' assumption.

\section{Schmidt  games} \label{games}

Wolfgang M. Schmidt introduced the games which now bear his name in \cite{schmidt}.   The simplified account which we are about to present  is more than adequate for the purposes of this paper.

Suppose that $0<\alpha<1$ and $0 <  \beta < 1$.   Consider the following game involving players A and B. First, B chooses a  closed interval $B_1 \subset \RR$. Next,  A chooses a closed interval $A_1$ contained in $B_1$ of length $ \alpha | B_1 | $.  Then,  B  chooses at will a closed interval $B_2$ contained in $A_1$ of length $ \beta | A_1 | $.  Alternating in this manner between the two players, generates  a nested sequence of closed intervals in $\RR$:
\[
B_1\supset A_1\supset B_2\supset A_2\supset\ldots \supset B_m\supset A_m\supset \ldots
\]
with lengths
\[
|B_m| \, = \,   (\alpha \, \beta )^{m-1}  |B_1|  \text{\quad and \quad} |A_m| \, =  \,   \alpha \,   |B_m|      \, .
\]

\noindent A subset  $S$ of $\RR$ is said to be {\em $(\alpha,\beta)$-winning} if  A can play in such a way that the unique point of intersection
$$
\bigcap_{m=1}^{\infty}   B_m   \,  = \, \bigcap_{m=1}^{\infty}   A_m
$$
lies in $S$, regardless of how  B plays.  The set $S$ is called {\em
$\alpha$-winning} if  it is  $(\alpha,\beta)$-winning for all
$\beta\in (0,1)$. Finally, $S$ is simply called {\em winning} if it
is $\alpha$-winning for some~$\alpha$.  Informally, player B tries
to stay away from the `target' set $S$ whilst player A tries to land
on $S$. The following results are due to Schmidt \cite{schmidt}.

\begin{SCH1}\label{sch_th1}
If $S \subset \RR$ is an
$\alpha$-winning set, then   $\dim S=1$.
\end{SCH1}

\begin{SCH2}\label{sch_th2}
The intersection of countably many  $\alpha$-winning sets is
$\alpha$-winning.
\end{SCH2}

Armed with these statements it is obvious that
$$
{\rm Theorem \ }  \ref{SVmixschdt} \quad  \Longrightarrow \quad {\rm Theorem \ } \  \ref{mixschdt}  \ .
$$

\section{Proof of Theorem~\ref{SVmixschdt}}

For any real $c >0$, let
$\bad_\DDD(c;i,j)$ be the set of $x\in \RR$  such that

\begin{equation}\label{neq2}
\max\{|q|_\DDD^{1/i}, ||qx||^{1/j}\} >  \frac{c}{q}   \qquad \forall
\ q \in\NN .
\end{equation}

\noindent It is easily seen that $\bad_\DDD(c;i,j)$ is a subset of $\bad_\DDD(i,j)$. Moreover,   it
has a   natural  geometric interpretation in terms of avoiding neighbourhoods of rational numbers.  Let
$$
\CCC_c :=  \{ r/q   \in \QQ  \, : \,   (r,q)=1  \, , q > 0  \ \  {\rm and  \ }  \     |q|_\DDD<c^i q^{-i}   \}
$$
and
\[
    \Delta_c(r/q):=\left[\frac{r}{q}-\frac{c^j}{ q^{1+j}
} \, , \, \frac{r}{q}+\frac{c^j}{ q^{1+j} }\right]
\]
Then,
$$
\bad_\DDD(c;i,j)  =   \{x \in \RR :  x \notin \Delta_c(r/q)  \ \  \forall \ r/q \in \CCC_c \}  \ .
$$

\vspace*{3ex}

Now with reference to \S\ref{games},  let $\bad_\DDD(i,j) $ be the target set
$S$  and $\alpha \in (0,1)$ be a fixed real number at our disposal. Suppose player
B has chosen some  $\beta \in (0,1) $  and an interval $B_1$.  Let
\[
    R:=\frac{1}{\alpha\beta} \, >  1 \,
\]
and fix  $ c > 0 $ such that

\begin{equation}\label{ineq_c}
c  \, < \,  \min \left\{ ( 4 \,  R \, |B_1|^{-1} )^{-1/j}    \,  ,
\,
 (  2  \, R^{\frac{i}{j+1}} \,   |B_1| )^{-1/i}  \, \right\}   \, .
\end{equation}

\noindent By definition, for each $m \geq 1$
\[
|B_m| \, = \,   R^{-m+1}  |B_1|     \, .
\]

\noindent The  `winning' strategy that  player  A  adopts is as follows.  {\em If $B_m$ is the interval player A inherits from player B, then A will
choose an interval $A_m \subset B_m$ with $ |A_m|  =     \alpha \,   |B_m|  $ such that }
\begin{equation}\label{cond_jn}
A_m\cap \Delta_c(r/q)=\emptyset\quad \forall\ r/q\in \CCC_c  \quad {with \ } \ 0 < q^{1+j}<R^{m-1}\; .
\end{equation}

\noindent Suppose for the moment that player A can adopt this
strategy with $\alpha = 1/4$. Then $$ \bigcap_{m=1}^{\infty}    A_m
\, \in \bad_\DDD(c;i,j)  \; \subset \; \bad_\DDD(i,j) \,     $$ and
it follows that  $\bad_\DDD(i,j)$   is $1/4$--winning as claimed.
We use induction to prove that such a strategy exists.

For $m =1$,  player $A$ can trivially  choose an interval $A_1$ satisfying \eqref{cond_jn}  since there are no rationals with $ 0<q < 1 $.    Now suppose the intervals
\[
B_1\supset A_1\supset B_2\supset A_2\supset\ldots \supset B_n\supset A_n\supset  B_{n+1}
\]
have been determined with each   $A_m$  ($1 \le m \le n$) satisfying  \eqref{cond_jn}. The goal is to show that there exists an interval  $A_{n+1}$  satisfying  \eqref{cond_jn}.   To begin with observe that since $B_{n+1} $ is nested in  $A_n$, we have that
\begin{equation*}\label{cond_in}
B_{n+1}\cap \Delta_c(r/q)=\emptyset\ \ \forall\; r/q\in \CCC_c\ \mbox{
with }\ 0< q^{1+j}<R^{n-1}.
\end{equation*}
Thus, since $ A_{n+1}$ is to be nested in $B_{n+1}$, it  follows that $ A_{n+1}$ will satisfy  \eqref{cond_jn} if
\begin{equation}\label{cond_jnsv}
A_{n+1}\cap \Delta_c(r/q)=\emptyset\quad \forall\ r/q\in \CCC_c(n) \ ,
\end{equation}
where
$$
\CCC_c(n):=\{r/q\in \CCC_c\,:\, R^{n-1}\le q^{1+j}<R^n\}.
$$

\vspace*{3ex}

\noindent{\em Fact 1. }  Let  $r/q\in \CCC_c(n)$. Then
\begin{eqnarray*}
|\Delta_c(r/q)|  & =  & \frac{2c^j}{q^{1+j}}
  \ \le \ 2c^j  \,   R^{-n+1} \\[1ex]
& \stackrel{\eqref{ineq_c}}{<} &      \, \textstyle{\frac12} \; |B_{n+1}|     \, .
\end{eqnarray*}
%The upshot is that there exits  an interval contained within  $ B_{n+1}$ of length
%$\frac{1}{4} |B_{n+1}| $ which avoids $\Delta_c(r/q)$.

\vspace*{3ex}

\noindent{\em Fact 2. } Let  $r_1/q_1,r_2/q_2\in
\CCC_c(n)$. Then there exists non-negative integers $k_1$ and $k_2$ such that
$$
q_s  \, =  \,  D_{k_s}q^*_s    \quad {\rm and }  \quad  \ q_s  \not \in D_{k_s+1}  \, \ZZ \qquad  (s=1,2) \, .
$$
Since $|q_s|_\DDD  \, <   \, c^i \, q_s^{-i}$, we have that
$$
D_{k_s}  \; >  \; c^{-i}  \,  q_s^i  \;  \ge \; c^{-i}  \, R^{\frac{(n-1)i}{j+1}} \, .
$$
Hence, it follows that
\begin{equation*}\label{ineq_gcd}
(q_1,q_2)> c^{-i}R^{\frac{(n-1)i}{j+1}}
\end{equation*}
and so
\begin{eqnarray*}
\left|\frac{r_1}{q_1}-\frac{r_2}{q_2}\right|  & \ge  &
\frac{(q_1,q_2)}{q_1q_2} \, >  \,
c^{-i} R^{\frac{(n-1)i}{j+1}} R^{- \frac{2n}{j+1}}   \\[1ex]
& = & c^{-i}  \, R^{-\frac{i}{j+1}} \, R^{-n}  \\[1ex]  &  \stackrel{\eqref{ineq_c}}{>} &      \, 2|B_{n+1}|     \, .
\end{eqnarray*}

\vspace*{3ex}

 A  straightforward consequence of the above two facts is that there is at most one rational  $r/q\in \CCC_c(n)$ such that
$$
\Delta(r/q)  \, \cap \, B_{n+1}  \neq \emptyset \,   .
$$
This together with Fact 1 implies that there is at least one interval $I \subset B_{n+1} $  of length $\frac{1}{4} |B_{n+1}| $  that  avoids $\Delta_c(r/q)$ for all $r/q\in \CCC_c(n)$.    With $\alpha = 1/4$, player A takes $A_{n+1}$ to be any such interval and this completes the induction step.  The upshot is that for $\alpha=1/4$ and any  $\beta \in (0,1)$
there exists a winning strategy for player A.

\vspace*{3ex}

\noindent{\em Remark. }
The arguments used to prove Theorem \ref{SVmixschdt}  can be naturally  modified to establish the following statement.  For any given finite number of  sequences  $\DDD_1,\ldots,\DDD_s$  and strictly positive real numbers $  i_1, \ldots,i_s, j $  satisfying $  i_1+\ldots+i_s+j=1 $,
the set of $x  \in \RR $  such  that
$$
\max\{  \; |q|_{\DDD_1}^{1/i_1}, |q|_{\DDD_2}^{1/i_2},\ldots,
|q|_{\DDD_s}^{1/i_s},||qx||^{1/j}  \, \}  \, >  \,  c(x) \ q^{-1} \quad \forall \  q \in \NN
$$
 is $1/4$--winning.

\section{The genuine mixed Schmidt conjecture?}

If $s=0$, then let us adopt  the convention that $ x^{1/s} :=0 $. Then $\bad_\DDD(0,1)$ is identified with  the standard set $\bad$  of badly approximable numbers   and $\bad_\DDD(1,0)$ is identified with  $\RR$. With this in mind, we are able to replace  `$< $'   by  `$\le$' in \eqref{neq1}  without effecting the statements of Theorems \ref{SVmixschdt} \& \ref{mixschdt}.   Moreover, it enables us to consider the following  generalization of Schmidt's conjecture.  Given    $\mathcal{D} $   and real numbers $i,j,k$ satisfying
\begin{equation}\label{neq11}
0\le i,j,k \le 1   \quad {\rm and  } \quad  i+j+k=1  \, ,
\end{equation}
let  $\bad_{\DDD}(i,j,k)$ denote the set
of   $(x,y) \in \RR^2$ for which  there exists a  constant $c(x,y)> 0$ such that
\begin{equation*}
  \max \{  \, |q|_{\DDD}^{1/i} \; ,  \; \|qx\|^{1/j} \; , \ \|qy\|^{1/k} \,  \} \, > \,
  c(x,y) \ q^{-1} \quad \forall \  q \in \NN   \ .
\end{equation*}

\noindent Naturally,   $\bad_{\DDD}(1,0,0) :=  \RR^2  $,  $\bad_{\DDD}(0,1,0) :=  \bad \times \RR  $,
$\bad_{\DDD}(0,0,1) :=  \RR \times \bad $,   $\bad_{\DDD}(0,j,k) :=  \bad(j,k)$, $\bad_{\DDD}(i,0,k) :=  \bad_{\DDD}(i,k)$ and  $\bad_{\DDD}(i,j,0) :=  \bad_{\DDD}(i,j)$.

\vspace*{1ex}

\begin{PMSconj}
 For any $(i_1,j_1,k_1)$ and
$(i_2,j_2,k_1)$ satisfying \eqref{neq11}, we have that
 $$ \bad_\DDD(i_1,j_1,k_1) \cap \bad_\DDD(i_2,j_2,k_2) \  \neq \  \emptyset \ .
 $$
\end{PMSconj}

Observe that when  $i_1=i_2=0$, this `mixed' conjecture reduces to the classical Schmidt conjecture. On the other hand, when $j_1=j_2=0$ or $k_1=k_2=0$ the above conjecture reduces to the mixed Schmidt conjecture investigated in this paper. In view of the results established  to date  it is reasonable to expect that the  following is true.

\vspace*{1ex}

\begin{MBPV}
Let  $(i_t,j_t,k_t)$ be a countable number of triples of real numbers
satisfying \eqref{neq1}.
Then
$$
\dim \Big(\bigcap_{t=1}^{\infty} \bad_{\DDD}(i_t,j_t,k_t) \Big) = 2   \ .
$$
\end{MBPV}

\noindent Note that this conjecture is open even when $i_t = 0 $ for all $t \in \NN$. The point is that this situation is not covered by
Theorem~BPV  since we have not imposed the  condition that  $ \liminf_{t \to \infty }  \min\{ j_t,k_t  \}  > 0   $.

\vspace{10mm}

{\small

\noindent Dzmitry A. Badziahin: Department of Mathematics,
University of York,

\vspace*{-1mm}

\noindent\phantom{Dzmitry A. Badziahin: }Heslington, York, YO10 5DD,
England.

%\vspace{0mm}

\noindent\phantom{Dzmitry A. Badziahin: }e-mail: db528@york.ac.uk

\vspace{5mm}

\noindent Jason Levesley: Department of Mathematics, University of
York,

\vspace*{-1mm}

\noindent\phantom{Jason Levesley: }Heslington, York, YO10 5DD,
England.

%\vspace{0mm}

\noindent\phantom{Jason Levesley: }e-mail: jl107@york.ac.uk

\vspace{5mm}

\noindent Sanju L. Velani: Department of Mathematics, University of
York,

\vspace*{-1mm}

\noindent\phantom{Sanju L. Velani: }Heslington, York, YO10 5DD,
England.

%\vspace{0mm}

\noindent\phantom{Sanju L. Velani: }e-mail: slv3@york.ac.uk

}

\end{document}